\documentclass[11pt,a4paper]{article}

\usepackage[utf8]{inputenc}
\usepackage{amsmath}
\usepackage{amsfonts}
\usepackage{amssymb}
\usepackage{makeidx}
\usepackage{graphicx}
\graphicspath{ {./} }
\usepackage{listings}
\usepackage{tocloft}
\usepackage{verbatim}
\usepackage{xcolor}

\usepackage{calc}
\usepackage{CJKutf8}
\usepackage[hidelinks]{hyperref}

\usepackage{epigraph}
\newtheorem{le[mm]a}[theorem]{Le[mm]a}

\newcommand{\qed}{\nobreak \ifvmode \relax \else
      \ifdim\lastskip<1.5em \hskip-\lastskip
      \hskip1.5em plus0em minus0.5em \fi \nobreak
      \vrule height0.75em width0.5em depth0.25em\fi}
      
       \usepackage[polutonikogreek,english]{babel}
    
\usepackage[normalem]{ulem}

\usepackage{tikz, tkz-tab, pgfkeys}
\usetikzlibrary{calc,decorations.pathmorphing,shapes}

  \newcounter{exer}[section]

    \makeindex

\title{André and  Simone Weil:  Mathematics, social activism and Indian culture}
\author{Athanase Papadopoulos and Susumu Tanabé}
 
\begin{document}

\maketitle

\begin{abstract}
This is an essay on the relation of André and Simone Weil  with Indian culture and Sanskrit literature, especially the \emph{Bhagavad G\=\i t\=a}, a  Hindu scripture which they knew well, which they quoted extensively, and which guided them in making important life decisions. In addressing this question, we will also talk about the life paths of the two Weils, and more specifically about certain aspects that relate to their deep convictions.

The final version of the article will appear in the Handbook of the Mathematics of the Arts and Sciences (New edition), ed. B. Sriraman, Springer,  2026.

\medskip


\medskip


\end{abstract}

\section{Introduction}

André Weil is \index{Weil, André} known in the mathematical community as a major figure of the twentieth-century. His sister, Simone Weil,\index{Weil, Simone} was a philosopher and mystic, and is much better known than he is in the general cultural world. In this article, we focus on André and Simone Weil's affinity with Indian culture and Sanskrit literature, with which they had an especially close relationship, in particular with the 
\emph{Bhagavad G\=\i t\=a}, which they knew well and quoted extensively, and for which they expressed a true devotion. We shall address the question of how reading this Hindu scripture was linked to their intimate convictions and to important decisions they made in their lives. 
Before we do so, however, it is necessary to properly situate André and Simone Weil, both in relation to each other and within the society in which their lives evolved.

 \section{André and Simone Weil}
André Weil\index{Weil, André} was born in Paris on May 6, 1906 and died in Princeton on August 6, 1998. Several written sources of information will perpetuate his memory. His own wide ranging writings contain an extensive record of his life and work and will help preserving them in the cultural memory. Among these writings we mention especially his reflections on mathematics and history of mathematics, the commentaries on his works that he wrote for his \emph{Oeuvres Scientifiques --- Collected Papers} published in two volumes by Springer \cite{Weil-Works}, his autobiography, \emph{Souvenirs d'apprentissage}, translated into English under the title \emph{The Apprenticeship of a Mathematician}  \cite{Weil-Souvenirs}, his correspondence with his family and especially with his sister, published by Gallimard in vol. 1 of Tome VII of the latter's \emph{Complete Works} \cite{SWeil-Corresp}, and his correspondence with mathematicians, in particular with Henri Cartan\index{Cartan, Henri} \cite{Cartan-Weil}. There are numerous tributes to him written by colleagues who had known him and were familiar with his work, including Borel  \cite{Borel1, Borel2, Borel3, Borel4}, Cartan  \cite{Cartan}, Cartier  \cite{Cartier}, Iyanaga  \cite{Iyanaga}, Serre  \cite{Serre},  Shimura \cite{Shimura}  and others (see also the reports by Pekonen \cite{Pekonen} and Varadarajan  \cite{Varadarajan}). Another important source of information on André Weil\index{Weil, André} is the biography of his  sister, Simone Weil,\index{Weil, Simone} written by her friend Simone Pétrement\footnote{Simone Pétrement (1907-1992) was a French philosopher and essayist who had been Simone Weil's classmate at the Lycée Henri IV and the \'Ecole Normale Supérieure. Pétrement was one of the first three women to enroll the Literature section of the \'Ecole Normale. It was to her that the Weil parents entrusted the task of writing their daughter's biography and that André Weil confided his recollections of his sister. The biography she wrote, \emph{La vie de Simone Weil} \cite{Petrement}, was awarded a prize by the Académie Fran\c caise} which gives an affectionate picture of the family and intellectual atmosphere in which the two Weils, André and Simone, grew up.
Finally, let us mention the novel by André Weil's daughter, Sylvie Weil, titled \emph{Chez les Weil --- André et Simone}, with an English translation titled \emph{At Home with André and Simone Weil} \cite{Sylvie-Weil}.  The novel is about Sylvie Weil's family, in particular her father and her aunt, Simone, who are our two main subjects in this article.

There is probably no need to introduce André Weil\index{Weil, André} to the mathematician reading this article. Let us only recall that he is the main founder of the Bourbaki group\footnote{Cartier, in \cite{Cartier}, writes: ``Bourbaki\index{Bourbaki, Nicolas} was undoubtedly Weil's creation, even if there was no shortage of strong personalities." [In this article, all the translations from the French are ours, except when an official translation is indicated in the bibliography.]} and that his name is attached to a number of objects and theorems that had a profound and durable influence on mathematics,  among which we mention the
  Weil--Petersson K\"ahler metric on Teichm\"uller space, 
  the de Rham--Weil theorem  in algebraic topology, the Chern--Weil homomorphism in the theory of characteristic classes, the Taniyama--Shimura--Weil conjecture in number theory (a conjecture which became later the Modularity theorem, of which the Fermat last theorem is a special case) and the Weil conjectures on zeta functions in the theory of algebraic varieties over finite fields, conjectures that have had a considerable influence on the foundations of modern algebraic geometry. There are several other mathematical concepts called after André Weil.
 
 In this article, we shall talk about André  and Simone Weil's\index{Weil, Simone} attachment to India and Indian culture. Since we started this article by talking about André Weil, let us mention here that he spent two  years in India, 1930-1932, teaching at Aligarh Muslim University.  It is interesting to know how he came to hold this position, since this is related to the heart of our subject here, and we shall say now a few words on this.

 André Weil\index{Weil, André} writes in his \emph{Apprenticeship} that at the time he was working for his doctorate, he was highly attracted by India and its culture. As a student at the \'Ecole Normale Supérieure, he had established a relationship with Sylvain\index{Levi@Lévi, Sylvain} Lévi,\footnote{\label{f:Levi} Sylvain Lévi  (1863-1935) was a French indologist and professor at  Collège de France (the highest possible position in French Academia).}  the famous Parisian Indologist. He was taking Sanskrit classes, and he had nurtured in himself the firm intention of traveling to India some day.  He writes in \cite[p. 58]{Weil-Souvenirs} that one day in 1929, Sied\index{Masood, Sied Ross} Ross Masood,\footnote{Sied Ross Masood bin Mahmood Khan (1889-1937) was the Vice-Chancellor of Aligarh Muslim University starting in 1929, that is, the year André Weil started working there. At the latter's arrival to Aligarth, he stayed in Masood's house for several weeks. Weil writes in his \emph{Apprenticeship}: ``At this time in India, it was beginning to be understood that an
exclusively English view of the world was not suitable for India. Indeed,
it was for precisely this reason that Masood had brought me to Aligarh."
Masood was knighted by the British Government in the 1933. He was linked to the British novelist E. M. Forster, who dedicated to him his novel \emph{A Passage to India} (1924).} the Vice-Chancellor  of Aligarh University came to Paris and informed Lévi,\index{Levi@Lévi, Sylvain} with whom he was in contact, about a vacancy of a professorship of French civilisation at that university. 
  Weil tells us that Lévi\index{Levi@Lévi, Sylvain} phoned him and asked: ``Are you serious about going to India?" Weil replied: ``Of course." Lévi then enquired whether Weil would accept to teach there French civilisation.  Weil replied: ``French or any other, I don't care; to go to India I'd teach anything they want." Eventually, the position opening was delayed, but in the meantime a chair of mathematics was declared vacant at the same university. Lévi\index{Levi@Lévi, Sylvain} and Weil immediately took advantage of the opportunity, and so it was that Weil set sail for India in early 1930, staying there for two years. 
  
Although Weil's\index{Weil, André} relationship with the University of Aligarh was conflictual (the university was in need of a complete overhaul and, as it turned out, Masood had invited Weil there precisely to put things in order), it is clear from what Weil writes in his \emph{Apprenticeship} that his stay in India was for him one of the most interesting, perhaps the most interesting, stay he made in his whole life. He writes \cite[p. 63]{Weil-Souvenirs}:
\begin{quote}\small
[\ldots] Thus, in January 1930, began a stay that was to last more than two
years, and which left me with a prodigious treasure trove of impressions.
This plethora of stimuli can be compared only to what a small child takes
in during the first years of life on this earth --- except that I was fully
cognizant from the start, and I still have vivid memories. In Sanskrit, the
brahmin is said to be ``twice born," \emph{dvija}: the second birth is conferred upon
him by the brahmanic cord. It was not simply as a joke that, shortly before
I was to return to France, my friend\index{Vijayaraghavan,Tirukkannapuram} Vijayaraghavan\footnote{Tirukkannapuram Vijayaraghavan (1902-1955) was a young mathematician 
who had studied with G. H. Hardy at Oxford in the mid 1920s and who was recruited at Aligarh University by Weil after his arrival there in 1930. As a matter of fact, in his \emph{Apprenticeship} \cite[p. 83]{Weil-Souvenirs}, Weil mentions the names of three young mathematicians he met during his stay in India: ``I had some regrets that my efforts had failed, but at least they had succeeded
in creating bonds of friendship among Vijayaraghavan, Kosambi,\index{Kosambi, Damodar Dharmananda} and
Chowla,\index{Chowla, Sarvadaman} the three young mathematicians who seemed to me the most
promising at the time."} girded me with such
a cord. Was it not the symbol of my second birth?
\end{quote}

Now to Simone Weil.\index{Weil, Simone}

While André Weil\index{Weil, André} is known in the mathematical community as one of the most influential mathematicians of the 20th century, in the non-mathematical cultured community the family name Weil evokes rather his sister, Simone Weil\index{Weil, Simone} (1909-1943), who is known as a socially committed intellectual, philosopher and historical essayist, activist for workers' rights, and mystic. She was a thinker who aimed to bring together the foundations of Ancient Greek religion, Hinduism,\index{Hinduism} and Christianity, and  the author of an \oe uvre that attempted to do so.

André Weil\index{Weil, André} and his sister shared the same ideas on several important issues. They admired and understood each other from childhood, and they maintained a rich correspondence when life later separated them.
They both possessed a formidable cultural background, including a profound historical vision and a thorough knowledge of ancient science and languages. Above all, they shared a demanding and uncompromising attitude in the respective societies and milieux in which they lived as well as an intransigent approach to the world of ideas in which they were evolving.
Reading their respective biographies and the letters they exchanged in adulthood, it clearly appears that they had a fused relationship, even when circumstances meant they ended up living in distant parts of the globe. 
In an interview with André Weil\index{Weil, André} concerning his relation with his sister, who was 2 years and 9 months younger than him, to the question ``Was there a strong intimacy between the two of you?", he responds: ``A very big one. As a child, my sister would spend her time imitating me.  My grandmother, who liked to speak German from time to time, said she was a \emph{Kopiermaschine}."\footnote{A copy machine.} (\cite[p. 10]{Sud}, also quoted in \cite[p. 25]{SWeil-Corresp}).  It was probably because he was still so much affected by the very early death of his sister (at age 34) that André Weil does not directly refer to her as much as would have been natural (or, at least, he claims not to mention her enough often)   in his autobiography written 48 years after her death, the \emph{Apprenticeship of a mathematician}. He writes in the \emph{Foreword} of that book (p. 11): ``My sister, too, is not much mentioned in these memoirs. Some time ago, in any case, I recounted my memories of her to Simone Pétrement,\index{Pétrement, Simone}
who recorded them in her fine biography, \emph{Simone Weil: A Life}. It would
be superfluous to repeat herein those many details on our childhood which
are included in Pétrement's book." In spite of this, it is easy to guess, while reading the \emph{Apprenticeship}, that Simone Weil\index{Weil, Simone} was always present in her brother's thought. She is also on most of this book's photographs.   
Let us quote another passage from the \emph{Souvenirs}, concerning his sister:

\begin{quote}\small
 As children, Simone and I were inseparable; but I was always the big brother and she the little sister. Later on we saw each other only rarely, speaking to one another most often in a humorous vein; she was naturally bright and full of mirth, as those who knew her have attested, and she retained her sense of humor even when the world had added a layer of inexorable sadness. In truth we had few serious conversations. But if the joys and sorrows of her adolescence were never known to me at all, if her behavior later on often struck me (and probably for good cause) as flying in the face of common sense, still we remained always close enough to one another so that nothing about her really came as a surprise to me --- with the sole exception of her death. This I did not expect, for I confess that I had thought her indestructible. It was not until quite late that I came to understand that her life had unfolded according to its own laws, and thus also did it end. I was little more than a distant observer of her trajectory.
\end{quote}

Sylvie Weil, André Weil's\index{Weil, André} daughter, in the novel we mentioned, describes how her father, in moments of absent-mindedness, would confuse her with his sister and call her Simone, and that he sometimes unintentionally made his daughter feel, like a reproach, that she was far from the latter's intellectual level.  She writes: ``Genius was two-headed. My father had a double, a feminine double, a dead double, a ghostly double. For, yes, in addition to being a saint, my aunt was a double of my father, whom she resembled like a twin. A terrifying double for me, since I looked so much like her. I looked like my father's double.'' Sylvie Weil also recounts that at times, she was ashamed of her aunt, like the time she was at a dinner party and felt obliged to divert attention away from her aunt, because of a feeling of embarrassment towards the latter, so as not to answer a question from one of the guests asking if she was related to the philosopher ``who was admired by Catholics and hated Jews."

%
%
%
%

\section{Simone Weil, Christian spirituality and social activism}

Simone Weil\index{Weil, Simone} grew up under the protection of her brother, who was brilliant in his studies: he passed his \emph{baccalauréat}\footnote{\label{f:Bac} The  \emph{baccalauréat} is the French national high-school diploma. André Weil,\index{Weil, André} passed the two baccalaureates, sciences and humanities, which he took in succession. He got the \emph{Mention très bien}, which is the highest rating, for both.} with honors at age sixteen, was ranked first at the \emph{agrégation} in mathematics at age nineteen,\footnote{\label{n:Agreg} The French \emph{agrégation} is a very selective French diploma entitling the holder to teach in a high school with substantial salary, fewer class hours and other privileges.    Most French mathematicians had started by passing the agrégation. Once one  has successfully passed this competitive examination,  unless one continues for a PhD, he or she is appointed to a specific town (depending on ranking, and also on the current need for teachers). André Weil obtained his agrégation in Mathematics in 1925, at age 19, which is very exceptional.} and was awarded a doctorate at age 22. During the First World War, their father was mobilized as a military doctor and André\index{Weil, André} did his schooling at home, studying  alone from books. He taught his five-year-old sister to read on the sly; the two children wanted to surprise their father: one day when he came home, Simone took his newspaper and started reading it.
 
Let us move on to what Simone says about her relationship with her brother. We start with an excerpt from a letter written to Joseph-Marie Perrin,\index{Perrin, Joseph-Marie} a Catholic priest who had become her spiritual mentor. 
The letter is dated May 15, 1942, and was later published
as part of her \emph{Spiritual autobiography} \cite{SW-Attente}, a collection of letters she wrote\index{Weil, Simone} to Father Perrin between January and May 1942.\footnote{In her letters to J.-M. Perrin, Simone Weil makes explicit her attachment to the Catholic Church, but her refusal to get baptised, not least because of her allegience and faithfulness  to other spiritual traditions, especially those of the ancient Greek philosophers and Hinduism, an allegience for which the official Catholic Church would condemn her  for heresy. She also writes, in the same letter to Fr. Perrin, ``[\ldots] The pure and simple resolution not to think at all about the question of my possible entry into the Church [\ldots] It is quite possible that one day I will suddenly feel the irresistible impulse to ask for baptism at once, and I will run to ask for it."} The letter shows first a young Simone\index{Weil, Simone} going through an existential crisis, stunned by her brother's achievements and finding herself excluded from his world of ideas, but who would later gain the innermost certainty, that she would keep for the rest of her life, that any human being can attain truth if the intention and the effort to know are strong enough. She writes:

\begin{quote}\small
At age fourteen, I fell into one of those bottomless despairs of adolescence, and seriously thought about dying, because of the mediocrity of my natural faculties. 
The extraordinary gifts of my brother, who had a childhood and youth comparable to Pascal's, forced me to be aware of this. I had no regrets about external successes, but about not being able to hope for any access to that transcendental realm where only truly great men can enter, and where truth dwells. I would rather die than live without it. 

After months of inner darkness, I suddenly and forever had the certainty that any human being, even if his natural faculties are almost nil, can enter that realm of truth reserved for genius, if only he desires the truth and makes a perpetual effort of attention to reach it. In this way, he too becomes a genius, even if, for lack of talent, this genius is not visible from the outside. 
Later, when my headaches caused my few faculties to become paralyzed, a paralysis that I very quickly assumed was probably definitive, this same certainty made me persevere for ten years in efforts for attention that had almost no hope of results.

\end{quote}

  Simone Weil\index{Weil, Simone} enrolled the \'Ecole Normale Supérieure,\footnote{Simone Weil entered the \'Ecole Normale Supérieure at age 19. Her brother had entered there at age 16.} not in mathematics as her brother, but in the humanities section. At the \'Ecole, she, together with some of her fellow students, took part in political actions, primarily through petitions, gathering money for a union's strike fund or for unemployment funds, collecting signatures for an appeal against military preparation that granted officer rank through shorter military service, and she was involved in many other causes. These activities are recounted in detail in \cite{Petrement}.

 In September 1930, after obtaining the \emph{agrégation},\footnote{See Footnote \ref{n:Agreg}. Jean-Paul Sartre\index{Sartre, Jean-Paul} and Simone de Beauvoir,\index{Beauvoir@de Beauvoir, Simone} who were of Simone Weil's generation  and who, like her, passed the Philosophy agrégation,  both obtained it in 1929 (at age 24 and 21 respectively). Sartre was ranked first and Beauvoir second at this national exam. In fact, this was the second year Sartre had taken this exam; the first time he had failed.} Simone Weil\index{Weil, Simone} was appointed philosophy professor at a high school in of Puy-en-Velay, a provincial town in the Central South of France. There, in those difficult pre-war years, she very soon became active in the trade union movement, supporting workers' strikes against unemployment and falling salaries.  Shortly after starting her work as a teacher, she decided to live on the minimum wage paid to the unemployed in the town where she was teaching, donating the remainder of her salary to workers in need. She became a union activist, contributing essays in syndicalist and revolutionary journals, and she joined anarchist and communist circles, publishing articles in the Marxist periodicals \emph{La Critique sociale}, \emph{La Révolution prolétarienne} and other left-wing publications. 
 She was following, with much anxiety, Hitler's rise to power in Germany at the time, while also being critical of Stalin's repression policies in the Soviet Union. She once had a stormy discussion with Leo Trotsky,\index{Trotsky, Leo} in her parents' apartment where she had invited him. Trotsky was at that time Stalin's most famous opponent. He was living in exile and passing some time in Paris. With his wife, child and bodyguards, he stayed in Simone Weil\index{Weil, Simone} parent's apartment, at her invitation.    Even if she was attracted by Trotsky's ideals, Simone Weil\index{Weil, Simone} was more critical than he was concerning the effects of the October 1917 revolution, especially with regard to the one party regime and the bureaucracy that had taken hold in Russia and was crushing the working class. The discussion between Simone Weil\index{Weil, Simone} and  Trotsky is reported on in several writings of her and about her, in particular the booklet \emph{Conversation with Trotsky} \cite{Weil-Trotsky}. See also her comments in \cite{Weil-Trotsky2} on the article \cite{Trotsky} by Trotsky.\footnote{It is reported that the discussion was stormy and that it focused in particular on the role of the German Communist Party in triggering the revolution. Trotsky, who was surprised by Simone Weil's attitude, asked her: ``If you feel this way, why are you receiving us? Are you the `Salvation Army'?" It is also said that Natalia Sedova, Trotsky's wife, stunned by this discussion, exclaimed: ``Look at this child standing up to Trotsky!" The whole group finished the evening in a local movie theatre where the film \emph{Okraïna} by the Soviet director Boris Barnet was being shown, a film based on Konstantin Finn's eponymous novel.}
 This encounter is also mentioned in \cite{Cartier} and \cite{Weil-Souvenirs}.\footnote{In \cite[p. 97]{Weil-Souvenirs}, André Weil,\index{Weil, André} writes: ``My sister, pursuing her eventful career, had
begun to receive her German friends fairly often at our parents' home. Our
parents were sometimes troubled by these lengthy visits. Most of these
Germans were dissident socialists or communists who had fled from Hitler.
Unbeknownst to me at the time, even Trotsky stayed in my studio once, at
the end of 1933."}

At some point, in Le Puy-en-Velay, Simone Weil\index{Weil, Simone} joined a group of unemployed workers at a town council meeting and she spoke on their behalf. For this and more generally for her publicly expressed strong political views, the school where she was teaching asked for her replacement.  At that time, André Weil,\index{Weil, André} was working in India, teaching at Aligarh University.\index{Aligarh Muslim University}
He relates the following episode in his \emph{Apprenticeship}
\cite[p. 93]{Weil-Souvenirs}: 

\begin{quote}\small
While I was in India, my mother had kept me up to date on the eventful
life led by my sister (the trollesse, or female troll, as we called her within
the family) in Le Puy, where she had obtained her first appointment as a
philosophy professor in a girls' \emph{lycée}. [\ldots] because Simone had been so bold as to shake hands with unemployed
workers in a public square, and then to accompany their delegation to the
town council, the school administration had threatened her with disciplinary
action. This news had reached me just when I was fighting my own
battles with the administration at Aligarh, and naturally I was enchanted
with the way my sister stood up to the authorities. I had written her a letter
of congratulations, saluting her as ``Amazing Phenomenon" --- to which she
replied addressing me philosophically as ``noumenon."\footnote{In the Kantian philosophical language (\emph{Critique of Pure Reason}), the term \emph{noumenon} designates an object of pure thought  which exists independently of our sensory perception, as opposed to  the \emph{phenomenon}, which is the object of our sensible experience. The two Weils used to play with words such as these to address each other.} 
\end{quote}
Indeed, the letter, written by André Weil\index{Weil, André} to his sister, dated Nov. 5, 1931, starts with: ``Amazing Phenomenon, 
I heard about your achievements and could not resist sending you my warmest congratulations" \cite[p. 524]{SWeil-Corresp}.  In another letter, dated January 14, 1932 \cite[p. 525]{SWeil-Corresp}, also sent from Aligarh, André\index{Weil, André}  Weil calls his sister\index{Weil, Simone}  ``Champignon sur l'humus'' (fungus on top of humus), in reference to an article in a conservative 
regional newspaper that claimed that ``intellectuals, who want to `make it', grow on the misery of the poor world like mushrooms on humus.''
He writes in that letter:

\begin{quote}\small
It has been, I think, twenty-three years since you entered the phenomenal world, much to the annoyance of rectors and headmistresses. But I wouldn't take up my pen on such a trivial occasion if I didn't have more worthy reasons to send you my warmest congratulations and encouragement to persevere in the path you've opened up for yourself. As for my hoaxes, which are so abundant, my letters from Paris will give you all the details you need. The University of Aligarh is in disarray, it is a mess, and I can at least flatter myself that I gave the signal for the debacle. I am strongly considering returning via China and Russia, and I am working on a Russian grammar. 
\end{quote}
He\index{Weil, André} concludes the letter by the words ``Inqilab Zindabad''.\footnote{The sentence  ``Inqilab Zindabad'', which, in Hindustani, means ``Long live the revolution", was uttered during the India freedom struggle in the 1920s.}

Simone replied, in a letter written at the end of January 1932 \cite[p. 426]{SWeil-Corresp}: ``I will never forgive you if you see the USSR before I do.''
 In the same letter, she writes:
``You are a disruptive element too in your own way, as I see --- you gently demolish the universities you pass through [\ldots]\footnote{André Weil had a number of conflicts and setbacks with the administration of the University of Aligarh which led to his ousting, but this is another story.} Here, miners' strike in preparation, but the situation is pathetic --- Depressed proletariat --- Increasingly violent polemics between union tendencies."
Like her brother, she ends her letter with ``Inqilab Zindabad'', with a drawing of a sickle and hammer above her signature. The letter is written on a Miners' Union letterhead paper.

In 1934, Simone Weil\index{Weil, Simone} resigned from her teaching position with the intention of working as a factory worker. The reason was that she wanted to feel physically the great fatigue of assembly-line work and  the oppression, social suffering and anguish of unemployment experienced by the workers. She first worked  at the Alsthom factory in Paris, a factory which in those times was a mechanical construction plant. After Alsthom, she worked as a laborer in the steel industry at the Forges of Basse-Indre on the outskirts of Paris, and later as a machinist at Renault. In each of the factories where she worked, she asked to be assigned to the most strenuous jobs. She wrote to a friend that she had been dreaming of doing this kind of work for ten years \cite[p. 299]{Petrement}. 
Let us quote her from a letter she\index{Weil, Simone} sent to Boris\index{Souvarine, Boris} Souvarine\footnote{Boris Souvarine (1895-1984) was a political activist, journalist, historian and essayist and a friend of Simone Weil. His real name was Boris Lifschitz (he chose the pseudonym Souvarine after a character in Emile Zola,\index{Zola, Emile}  a Russian immigrant and anarchist working in a coal mine in Northern France). Souvarine was initially a Communist activist, before being excluded from the PCF in 1924, after he criticized the leadership of the Communist Party in Russia. An outspoken critic of Stalinism, Souvarine is the author of an early biography of Stalin \cite{Souvarine}.} in 1935, in which she describes her working conditions: ``Dear Boris, I am forcing myself to write you a few lines, because otherwise I wouldn't have the courage to leave a written record of the first impressions of my new experience. The so-called friendly little company turned out to be, on contact, first a rather large company, and then above all a dirty, dirty company. In this dirty company, there is a particularly disgusting workshop: it's mine." \cite[p. 29]{Weil-Condition}

Needless to say, Simone's parents, who lived a comfortable life in Paris (her father was a doctor), were very worried about her. All this is recounted in detail in the biography written by Pétrement \cite{Petrement} and in other biographies, but most of all, in her diary, notes and correspondence, published in her \emph{Collected Works}.  
In addition, for several consecutive summers, she worked in farms, harvesting grapes, potatoes and other earth produce at a frenetic pace, always wanting to do more. She also worked with fishermen at sea. From a certain point on, she was immersed in Christian religious texts. Later, she would write to Father Perrin  
 \cite{SW-Attente}: ``[\ldots] ``Peasants and workers possess that incomparably sweet closeness to God which lies at the bottom of poverty, lack of social consideration and long, slow suffering."
 
Meanwhile, the Spanish Civil war\index{Spanish Civil war} had broken out, a war between the so-called Republicans, comprising the left, the far left and the anarchists, and the nationalists, comprising the right-wing parties and the fascists, also known as the Francoists. Simone Weil\index{Weil, Simone} went to Spain and joined the anarchists with the objective of fighting the fascists, after having persuaded her parents that she was going to Spain as a journalist. She was only able to stay for two months, and was repatriated to France after a domestic accident in which her foot was badly burned.\footnote{It seems that the accident happened in a kitchen, where she had been assigned. In fact, during her stay in Spain, she was refused assignments at the front because she was extremely clumsy and totally impractical; everyone knew she was not able to handle a rifle. The same thing happened again  when she asked to be sent to the front, in France, to fight with the Resistance against the Germans and the Vichy regime; her request was refused. }

  \section{Greek, Sanskrit and Bhagavad-G\=\i t\=a}
Simone Weil\index{Weil, Simone} was deeply attached to the cultures of Ancient Greece and of India. Readers will have gathered that it was her brother who first exhibited this attachment. We have already mentioned that during the First World War, André Weil\index{Weil, André} did most of his schooling at home. We learn from his \emph{Apprenticeship} that by 1917 (he was 11), he knew enough ancient Greek and mathematics, which he had learned on his own. 
He writes (p. 40):  ``When
school reopened in October of 1917, my family was in Laval.  I had taught myself enough Greek, and I already knew enough mathematics, to be admitted into the classical section of the third form.\footnote{He means the fourth last year before finishing high school.}  I no longer remember
whether it was there that I read the first book of the Iliad, or whether I read it by myself, for my own pleasure; in any case, that was how I discovered
poetry, and that it is untranslatable''. Later, his sister also became a fervent devotee of ancient Greek literature. Readings notes, translations from ancient Greek and fragments of studies on ancient Greek thought that she wrote were published in 1953 in a volume titled \emph{La source grecque} (The Greek source).
But we shall talk now about Sanskrit.

 It was during the same period that André Weil learned ancient Greek that he became interested in Sanskrit. We read, on p. 24 ff. of the \emph{Apprenticeship}:
\begin{quote}\small
I had at the time, and indeed I still have Annandale's\footnote{Charles Annandale\index{Annandale, Charles} (1843-1915) was a Scottish editor, especially of reference books. André Weil refers probably to his 
\emph{Concise imperial dictionary:
literary, scientific and technical, with pronouncing lists of proper names, foreign words and phrases, key to names in mythology and fiction, and other valuable appendices}, published around 1915.} English
dictionary, which includes an introduction to Indo-European linguistics and Grimm's\index{Grimm's law} Law\footnote{Grimm's Law describes a series of phonetic transformations that trace the evolution of the occlusive consonants of Proto-Germanic (the ancestor of today's Germanic languages) from those of the common Indo-European, during the 1st millennium B.C.}  as well as fairly detailed etymological information going as far back as Sanskrit. I dreamed of one day being able to read, in the original, the epic poems written in all these languages. My romantic notion
of these epics later led me to seek out Sylvain Lévi's advice.\footnote{See Footnote \ref{f:Levi}.} 
[\ldots]  Fortunately, too, my courses at the \emph{Classes préparatoires}\footnote{In the French system, these are the two years of study that follow the  baccalauréat (see Footnote \ref{f:Bac}) and prepare students for the entrance examinations to the ``grandes \'Ecoles", that is, the universities that have the rank of ``Superior Schools", which, in the case of André Weil  and his sister, is the \'Ecole Normale Supérieure. Incidentally, André Weil entered the \'Ecole Normale after one year of \emph{Classes préparatoires}, instead of two years, which is extremely rare.} did not fill all my
time that year. I began to study\index{Jordan, Camille} Jordan.\footnote{This is the mathematician Camille Jordan (1838-1922). In his \emph{Apprenticeship} \cite{Weil-Souvenirs}, André Weil recounts that the Lycée Saint-Louis, the year he was preparing the baccalauréat (see Footnote \ref{f:Bac}), awarded an endowment prize to the best student in the mathematics section. The prize went to him and he received, as a reward,  the three-volume \emph{Cours d'Analyse} by Jordan and the 
\emph{Treatise of Natural Philosophy} by William\index{Thomson, William (Lord Kelvin)} Thomson (later Lord Kelvin) and Peter Guthrie Tait.\index{Tait, Peter Guthrie} He writes then: 
``Thanks to Hadamard, I learned analysis from Jordan (infinitely better than learning it from Goursat, as most of my classmates did) and I was initiated into differential geometry
by Thomson and Tait. Naturally I did not read these hefty volumes immediately; I think, however, that I did begin with Jordan the following year."} Also, my precocious and romantic
attraction to Sanskrit gave one of my father's friends the idea of introducing
me to the leading scholar in the field of Indian studies, Sylvain Lévi.
[\ldots] When Sylvain Lévi received me at his home in  Rue Guy-de-La-Brosse, he said to me: ``There are three reasons for studying Sanskrit,"
and he enumerated them: I believe they were the Veda, grammar, and
Buddhism. ``Which of these is yours?" I didn't dare tell him that I was
impelled by none of the three, but simply by the naive notion I had
of Indian epic poetry. I had already learned --- where, I do not recall
--- the alphabet and one or two declensions, and I asked his advice on
how to proceed. He told me that the best textbook by far was\index{Bergaine, Abel-Henri-Joseph} Bergaigne's,\footnote{Abel-Henri-Joseph Bergaine (1838--1888) is a French Indologist and  Vedic philologist. He was professor of Sanskrit and Comparative Indo-European at the Sorbonne.} but it was out of print; otherwise I should obtain a copy of
Victor\index{Henry, Victor} Henry's\footnote{Victor Henry (1850-1907) is a French linguist who  held the chair of Sanskrit and Comparative Grammar at the Paris Faculty of Letters (Sorbonne).} textbook. This I bought immediately, and I lost no time in putting it to use.

\end{quote}

We are now in 1922, and André Weil\index{Weil, André} is a student at the \'Ecole Normale Supérieure, in Paris. He writes in the \emph{Apprenticeship} \cite[p. 40 ]{Weil-Souvenirs}:
 
 \begin{quote}\small
At the end of the
year, wishing to devote part of my vacation to reading a Sanskrit text, I went to Sylvain Lévi\index{Levi@Lévi, Sylvain} for advice. From a shelf in his library, he pulled a small volume bound in red velvet. It was a ``native'' (to use the term in vogue at the time) edition of the  \emph{Bhagavad G\=\i t\=a}. ``Read this,'' he told me. ``First of all, you cannot understand anything about India if you haven't read it'' --- here he paused, and his face lit up --- ``and besides,'' he added, ``it is beautiful.'' With the help of a dictionary lent to me by Jules\index{Bloch, Jules} Bloch,\footnote{\label{f:Bloch} Jules Bloch (1880-1953)  was a French Indologist. At the time André Weil followed his courses, he was teaching, at the \'Ecole Nationale des Langues Orientales Vivantes and at the 
\'Ecole Pratique des Hautes \'Etudes. He became later Sylvain Lévi's successor at the Collège de France, teaching there Sanskrit language and literature.} and of the English translation of \emph{Sacred Books of the East},\footnote{The \emph{Sacred Books of the East}, in  50-volumes, contains English translations of sacred texts of Hinduism, Buddhism, Taoism, Confucianism, Zoroastrianism, Jainism, and Islam. It was edited by the philologist and Orientalist Max Müller and published by the Oxford University Press between 1879 and 1910.} borrowed from the
\emph{\'Ecole} library, I read the \emph{G\=\i t\=a} from cover to cover [\ldots]. The beauty of the poem affected me instantly, from the very first
line. As for the thought that inspired it, I felt I found in it the only form of
religious thought that could satisfy my mind. My sister and I had been
brought up without any semblance of religious education or religious
observance. [\ldots]
Later on, I felt at home in India
thanks to the \emph{G\=\i t\=a}; and much later, it was also thanks to the \emph{G\=\i t\=a} that my
sister's way of thought did not seem at all foreign to me.
\end{quote}

During his second year at the \'Ecole Normale, André Weil\index{Weil, André} regularly attended\index{Bloch, Jules} Jules Bloch's\footnote{see Footnote \ref{f:Bloch}.} course on the \emph{Veda} at the \'Ecole des Hautes \'Etudes,\index{Meillet, Antoine} Meillet's\footnote{\label{f:Meillet} Antoine Meillet (1866-1936) was a specialist of Sanskrit, Romance languages, Slavic languages, Irish, Iranian and Armenian. Between 1902 and 1906, he was professor of Armenian at the \'Ecole des Langues Orientales, and between 1906 and 1936, professor of comparative grammar at the Collège de France.} lectures on Indo-European linguistics, and
Sylvain\index{Levi@Lévi, Sylvain} Lévi's\footnote{André Weil, while giving a few details about the teaching approach of Meillet and Lévi, writes: ``These two teachers were without
peer" \cite[p. 42]{Weil-Souvenirs}.} on K\=alid\=asa's  
\emph{Meghad\=uta}\footnote{K\=alid\=asa is a poet and playwright who wrote in Sanskrit between the 4th and 5th centuries. He is one of the most famous classical authors of Sanskrit literature. The Meghad\=uta\index{Meghaduta@\emph{Meghad\=uta}} is an important Sanskrit lyric poem by K\=alid\=asa.\index{Kalidasa@K\=alid\=asa} André Weil, in his \emph{Apprentceship} (p. 42), writes that ``this poem consists
principally of a long speech made by a Yaksha in love to the `messenger
cloud' (the Meghad\=uta of the title) which will relay it to his distant best-beloved."} \cite[p. 42]{Weil-Souvenirs}. Talking about the latter, Weil recalls, in \cite[p. 42]{Weil-Souvenirs}: ``Meillet, already nearly blind --- this was said to be from deciphering
numerous undecipherable texts from Central Asia --- would hold forth
extemporaneously, treating his audience to observations which even my
ignorant mind could tell were strikingly original. Sylvain Lévi\index{Levi@Lévi, Sylvain} read and
explained the K\=alid\=asa poem, verse by verse, in his slightly muffled voice."

About fifty years later, at the time he was composing his \emph{Apprenticeship}, André Weil\index{Weil, André} writes  \cite[p. 149]{Weil-Souvenirs}: ``And I continue, with as much enjoyment as ever, my reading of the
\emph{G\=\i t\=a}. I had often reread certain passages, but I had read it straight through from beginning to end only once, in 1923 or 1924. I am pleased with myself for remembering my Sanskrit well enough to read such a text."

\section{André Weil and his \emph{dharma}}

In this section,\index{Weil, André} we shall be talking about  
\emph{dharma},\index{dharma@\emph{dharma}}  a complex and multivalent word in Indian religions and spirituality. The few words that follow are aimed primarily at the Western reader unfamiliar with Indian culture.

The word\index{dharma@\emph{dharma}}  \emph{dharma} can designate a set of norms, which may be social, family or personal, or a natural and cosmic law which may be universal. It can also refer simply to a personal truth or reality. It also means one's social or moral responsibilities or duties. 
We quote, and we shall quote again in its context, a sentence by André Weil:\index{Weil, André} `` Gauguin's \emph{dharma} was painting. Mine, as I saw it in 1938, seemed
clear to me: it was to devote myself to mathematics as much as I was able."
 In any case, it is probably in the various uses of this term by André and Simone Weil,\index{Weil, Simone} in the lines that follow, that its meaning in the present context will become clearer.

In\index{dharma@\emph{dharma}} the \emph{Apprenticeship} \cite[p. 126 ff.]{Weil-Souvenirs}, André Weil\index{Weil, André} shares his thoughts on compulsory military service, especially in times of war. He tells us that, during his years at the \'Ecole Normale, he was already struck by the fact that in France, during the First World War, an entire generation of researchers was sacrificed, since even those who returned after four years of war were no longer capable of exercising any scientific research profession.\footnote{Weil writes in his \emph{Apprenticeship}: ``In 1914, the Germans had wisely sought to spare the cream
of their young scientific elite and, to a large extent, these people had been
sheltered. In France a misguided notion of equality in the face of sacrifice
--- no doubt praiseworthy in intent --- had led to the opposite policy, whose
disastrous consequences can be read, for example, on the monument to the
dead of the \'Ecole Normale."} He writes: ``This was a fate that I thought it my duty, or rather my \emph{dharma},\index{dharma@\emph{dharma}} to avoid. [\ldots] I was not a conscientious objector in the sense usually given to this word, that is, someone who does not believe in killing or even using instruments of death, whether he sees this as a universally valid precept or as his own personal  \emph{dharma}."  He makes it clear that this is not a question of pacifism, and that he feels himself to be as much far from an unconditional pacifist as from uncompromising patriots. He tells us that this course of action, in 1940, got him into a series of troubles, which he says caused his sister a crisis of remorse, as she considered that, because at some point she had pacifistic inclinations, she might have influenced him in his decision not to submit to the regulation of national service for all.
 He says that if there was one thing that influenced him, it was the example of the mathematician Carl Ludwig Siegel,\index{Siegel, Carl Ludwig} who had told him how he had once deserted in 1918, while serving as a soldier in the German army in Alsace. It seems that
Siegel told him: ``Dieser Krieg war nicht mein Krieg," which means: ``This war was not my war".
This passage on war ends with the words: ``In this light, I felt justified in thinking that by exempting myself
from military law, I was, to the small extent permitted me by circumstances,
taking charge of my own destiny."

 At the moment he was called up as a soldier, it was the \emph{G\=\i t\=a}\index{Bhagavad@\emph{Bhagavad G\=\i t\=a}} that helped André Weil\index{Weil, André} make a decision. In the  \emph{Apprenticeship} (p. 24), he writes: ``I had been deeply marked by Indian thought and
 spirit of the \emph{G\=\i t\=a},\index{Bhagavad@\emph{Bhagavad G\=\i t\=a}}  such as I felt capable of interpreting it.''
He is referring here to that passage of the \emph{G\i t\=a}, known to the Hindus, in which Arjuna, the very famed warrior, crippled with emotions, expresses to Krishna doubts whether or not to go to battle, because he does not understand the reason why he has to fight and to kill people from the other clan, who are also members of his family. This  was during the Mah\=abh\=arata  war\index{Mahabharata@Mah\=abh\=arata  war}, which pitted the Kauravas against the Pandavas, two groups from the same family. Krishna's\index{Krishna} response to Arjuna  can be summed up by the fact that the latter is a warrior, and that his duty is to fight. The question, thus, indeed concerns the faithfulness to one’s \emph{dharma}.

About his own \emph{dharma},\index{dharma@\emph{dharma}} André Weil\index{Weil, André} writes, in the \emph{Apprenticeship}, after again quoting the \emph{G\i t\=a}, that ``the only recourse is
for each one of us to determine as best he can his \emph{dharma}, which is his
alone. Gauguin's \emph{dharma} was painting. Mine, as I saw it in 1938, seemed
clear to me: it was to devote myself to mathematics as much as I was able.
The sin would have been to let myself be diverted from it"  \cite[p. 126]{Weil-Souvenirs}.

Weil's refusal to serve in the army was not properly interpreted  by some of his colleagues. Cartier would later write \cite{Cartier}: ``Some translated this crudely as `my brain is too precious to be exposed to bullets', or `to go to war would waste my invaluable time, while I am full of exciting discoveries'.  This was seen as nothing more than the egotistical translation, in a time of crisis, of the `pure mathematician's motto', which is quite common in our milieu, and which often serves as an alibi for shirking responsibility.” 
This had consequences on Weil's career, which we shall address later in this article.

In any case, because he refused to join the army, André Weil\index{Weil, André} spent five months in prison in Rouen, a city in Normandy. We might add in passing that he would later say that, from a mathematical point of view, these five months constituted one of the most productive periods of his life.  Sanskrit and the \emph{G\=\i t\=a}\index{Bhagavad@\emph{Bhagavad G\=\i t\=a}}  accompanied him\index{Weil, André} during the period of his incarceration.  Let us quote a few passages from letters he sent to his wife, \'Eveline, from prison.  These excerpts are reproduced in \cite{Weil-Souvenirs}. From a letter sent on March 30, 1940 \cite[p. 146]{Weil-Souvenirs}:
``[\ldots] And then there are always my Sanskrit books. I am reading the
\emph{G\=\i t\=a},\index{Bhagavad@\emph{Bhagavad G\=\i t\=a}}  in small doses as one ought to read this book. The more detail one
absorbs, the more one admires it."  On April 7 of the same year, he writes \cite[p. 147]{Weil-Souvenirs}:
``Here are some lines from
the \emph{G\=\i t\=a}\index{Bhagavad@\emph{Bhagavad G\=\i t\=a}}  that I like very much: `A leaf, a flower, a fruit, some water,
whoever dedicates it with love, this love offering I accept with the devotion
of his soul.'\footnote{\emph{Bhagavad G\=\i t\=a},  IX.26.} The God Krishna\index{Krishna} is speaking. In Tamil it is written of him, `The
bread we eat, the water we drink, the betelnut we chew, all this is our
Krishna.' It is almost impossible to translate all of this --- in everything to
do with gods our language is slanted toward the idea of a personal god, that
is totally foreign to the Indians [\ldots]"

The question would come up again later, when André Weil\index{Weil, André} had left France. His sister, who was part of the \emph{France Combattante},\footnote{This was the name given in 1942 to the military and political resistance movement founded in London by General de Gaulle against the German occupation of France.} in a letter sent on April 23, 1943 \cite[p. 516]{SWeil-Corresp}, asks him to join this organization, specifying that ``membership implies only the affirmation that it was right and proper, in June 1940, to proclaim France's continued participation in the war; something I have never doubted". André Weil\index{Weil, André} refused. He replied that ``by joining the F. C., one necessarily places oneself at the disposal of the military authority --- something I am loathe to do, and which would only be excusable, in my own eyes, if I decided to become a military man in the strict sense of the word."

Again, it was a question of \emph{dharma}, and André Weil\index{Weil, André}  quoted the \emph{G\=\i t\=a}:\index{Bhagavad@\emph{Bhagavad G\=\i t\=a}} Either you completely become  a soldier, or you remain, as you are now, a brahman, and in the latter case you devote your life to intellectual activities: teaching, learning and meditation --- you cannot be both at the same time. He writes in the \emph{Apprenticeship} \cite[p. 124]{Weil-Souvenirs}:
\begin{quote}\small
 It seems fitting for me to indicate here the motives behind my
resolution not to serve, though I fear this may turn into a long and confused
explanation. I thought I was being perfectly lucid at the time; but is one
ever perfectly lucid when making a decision of grave impact? I have never
believed in the categorical imperative. The Kantian ethic, or what passes
for it today, has always seemed to me to be the height of arrogance and
folly. Claiming always to behave according to the precepts of universal
maxims is either totally inept or totally hypocritical; one can always find
a maxim to justify whatever behavior one chooses. I could not count the
times (for example, when I tell people I never vote in elections) that I have
heard the objection: ``But if everyone were to behave like you ... " --- to which
I usually reply that this possibility seems to me so implausible that I do not
feel obligated to take it into account.

On the other hand, I had been deeply marked by Indian thought and
by the spirit of the \emph{G\=\i t\=a},\index{Bhagavad@\emph{Bhagavad G\=\i t\=a}} such as I felt capable of interpreting it. The law is
not: ``Thou shalt not kill", a precept which Judaism and Christianity have
inscribed --- to what avail? --- in their commandments. The G\=\i t\=a\index{Bhagavad@\emph{Bhagavad G\=\i t\=a}} begins with
Arjuna, ``filled with the deepest compassion," stopping his chariot between
two armies, and ends with his lucid acceptance of\index{Krishna} Krishna's injunction to
go to combat unflinchingly. Like everything else in this world, combat is
an illusion: he does not kill, nor is he killed, whoever has known the Self.\footnote{Bhagavad G\=\i t\=a, II: 19.}
In the absence of any universal recipe to prescribe everyone's behavior, the
individual carries within him his own \emph{dharma}.\index{dharma@\emph{dharma}} In the ideal society of the
mythical times of the Mah\=abh\=arata,\index{Mahabharata@Mah\=abh\=arata  war} the \emph{dharma} comes from the individual's
caste. Arjuna belongs to a caste of warriors, so his \emph{dharma} is to
go to combat. Krishna\index{Krishna} is the exceptional being, the divinity incarnate.
``Whenever \emph{dharma}\index{dharma@\emph{dharma}} declines and its opposite triumphs, then I reincarnate
myself,"\footnote{Bhagavad G\=\i t\=a,  IV.7.} he says in a famous verse which was once applied to Gandhi.
Krishna exists outside the \emph{dharma}.

\end{quote}

 The question for André Weil,\index{Weil, André} if he solved it for himself, wasn't as easy as that. Indeed, he had advised his sister, who was a fighter for social justice, to leave France too, which she did in 1942, but something she would later regret. She wrote him from England, in a letter sent on April 17, 1943 (he was in the USA): 
``I am torn more and more cruelly day after day by regret and remorse for having been weak enough to have followed your advice a year ago".\footnote{i.e. for having left France.} And she adds: ``As for you, if you were now in conditions favorable to mathematical work, I would certainly advise you to think only of mathematics, and that, definitively, if possible, until death."

 Meanwhile, André Weil,\index{Weil, André} was being criticized in France, by certain colleagues, for deserting. Cartier writes in \cite{Cartier} that Weil would have liked to settle to France, but was prevented from doing so: ``The `patriot' party", he writes, ``never gave up on him; 
[\ldots]
 The most relentless member of the patriot party was Jean Leray.
[\ldots] 
  Leray worked hard to keep Weil out of the Sorbonne and the Collège de France. Weil stayed in Chicago, and my generation lost a master!"

To conclude this section, we quote an  excerpt from Antoine de Saint-Exupéry's\index{Saint-Exupéry@de Saint-Exupéry, Antoine} \emph{Terre des Hommes},\footnote{The title of the  English translation is \emph{Wind, Sand and Stars}.} in which the author speaks about Truth. The ideas in this passage ares very close in substance to what André Weil,\index{Weil, André} writes about his\index{dharma@\emph{dharma}} \emph{dharma}:\footnote{We are thankful to Alexey Sossinsky for pointing out to us this passage. Antoine de Saint-Exupéry (1900-1944) was one of the French intellectuals who, like André Weil, left France for the United States during the Second World War.  Among the others who, like him, left for the US, we mention André Breton, Claude Lévi-Strauss, Ossip Zadkine, Marc Chagall, Jean Renoir, and the mathematicians Jacques Hadamard and Claude Chevalley (the latter was already at Princeton when the war broke out), and there are many others.} ``Truth is not that which can be demonstrated by the 
aid of logic. If orange-trees are hardy and rich in fruit in 
this bit of soil and not that, then this bit of soil is what 
is truth for orange-trees. If a particular religion, or culture, or scale of values, if one form of activity rather 
than another, brings self-fulfillment to a man, releases the 
prince asleep within him unknown to himself, then that 
scale of values, that culture, that form of activity, constitute his truth. Logic, you say? Let logic wangle its 
own explanation of life."

\section{In guise of a conclusion}

 Throughout her life, Simone Weil,\index{Weil, Simone} like her brother, refused to make concessions. J.-M. Perrin writes that ``we cannot fully understand her if we forget the sharp character that emerged from her thirst for the absolute." He quotes her saying that ``the greatest evil is not evil but the mixture of good and evil" \cite[p. 29]{Perrin}. 
  To the same, Simone Weil writes, in a letter dated May 15, 1942 \cite{SW-Attente}: ``In the spring of 1940, I read the \emph{G\=\i t\=a}.\index{Bhagavad@\emph{Bhagavad G\=\i t\=a}} Strangely enough, it was while reading these marvellous words, so Christian-sounding, put into the mouth of an incarnation of God, that I felt strongly that we owe religious truth something more than the adherence to a beautiful poem, a kind of adherence that is much more categorical."  On June 9, 1943, that is, a month and a half before her death, in a letter to her parents sent from London where she had been exiled against her will, Simone Weil\index{Weil, Simone} writes: 
   ``I started doing a few lines of Sanskrit every day, in the G\=\i t\=a.\index{Bhagavad@\emph{Bhagavad G\=\i t\=a}} It feels so good to be in contact with Krishna's\index{Krishna} language." At the end of the same letter, she tells her mother: ``Enjoy the beautiful days and think of Krishna".\index{Krishna} Until the last moment, she was hiding the fact that she was ill and with an extreme form of anorexia.
She passed away in a sanatorium in Ashford (England), 
  on August 24, 1943, aged 34. She had contracted tuberculosis, but what led her to death was physical and moral exhaustion, and starvation: she could not accept to get nourished normally while her compatriots were starving under Nazi occupation.  
   Almost all her works were published posthumously. From 1949 onwards, Albert Camus\index{Camus, Albert} was responsible for publishing some of her most important writings. In a letter to her parents, he wrote that their daughter is ``the only great spirit of our times."  The complete works of Simone Weil are published in seven tomes distributed over sixteen large volumes.

After Chicago, André Weil\index{Weil, André} joined Princeton's Institute for Advanced Study, where one of his colleagues, J. Robert Oppenheimer,\index{Oppenheimer, Julius Robert}\footnote{J. Robert Oppenheimer was professor at the Institute for Advanced Study for the period 1947 until 1967. He served as Director of the Institute from 1947 until 1966, thus far the longest tenure of any Institute Director, and it was at a time when André Weil was a member.} was also\index{Weil, André} a lover of the Baghavad-G\=\i t\=a,\index{Bhagavad@\emph{Bhagavad G\=\i t\=a}} which he also read in Sanskrit. We know the famous verse ``Now I am become Death, the destroyer of worlds'', and the fatal moment of the test known as the Trinity test,\index{Trinity test} July 16, 1945, when Oppenheimer recited it, but that is a separate subject.

 Let us conclude by quoting André Weil,\index{Weil, André} from his memories of his prison in Rouen \cite[p.142]{Weil-Souvenirs}:
``[\ldots] The door, like that of a safe but pierced with the usual peek-hole, opened with a loud rattling of keys and a grinding of hinges --- aptly illustrating the expression `as graceful as a prison gate.' I had with me the \emph{Bhagavad G\=\i t\=a},\index{Bhagavad@\emph{Bhagavad G\=\i t\=a}} the \emph{Chandogya-Upanishad},\index{Chandogya-Upanishad} novels by Balzac,\index{Balzac@de Balzac, Honoré} and, on my sister's enthusiastic recommendation, Retz's memoirs.\footnote{Cardinal de Retz,\index{Retz@de Retz, Cardinal (Jean François Paul de Gondi)} 
 Jean François Paul de Gondi (1613-1679) is a French ecclesiastic, known as a writer for his \emph{Mémoires}, an invaluable document for the historical and sociological information it contains.} I also had what I needed to resume my mathematical work.
[\ldots] I often recited to myself the lines of the  \emph{G\=\i t\=a}:\index{Bhagavad@\emph{Bhagavad G\=\i t\=a}}  `Patram
puspam phalam toyam \ldots' (A leaf, a flower, a fruit, water, for a pure heart everything can be an offering)."\footnote{\emph{Bhagavad G\=\i t\=a}, IX.26.}
 
 \medskip

\noindent{\it Acknowledgement.} The first author gratefully acknowledges invitations to several institutions in India where this paper was written: Banaras Hindu University, Jawaharlal Nehru University in Delhi, the Tata Institute for Fundamental research in Bombay and the International Centre for Theoretical Sciences of Bangalore. Discussions with Shubhabrata Das were especially enlightening. We  thank Arkady Plotnitsky, Alexey Sossinsky and Bharath Sriraman who took the trouble to read carefully a preliminary version of this text, to give their impressions and to ask for clarifications.  We also thank Himalaya Senapati for his help in Sanskrit.

\bigskip

\noindent {\bf Authors' addresses:}

\noindent Susumu Tanabé,
Moscow Institute of Physics and Technology,
Department of Discrete Mathematics,
Dolgoprudny, Moscow Region, Russian Federation.

\noindent  email : tanabe.s@mipt.ru

\medskip

\noindent Athanase Papadopoulos,  Université de Strasbourg et  Centre National de la Recherche Scientifique,
Institut de Recherche Mathématique Avancée,
7 rue René Descartes,
67084 Strasbourg Cedex France;

 \noindent Centre de Recherche et d'Expérimentation sur l'Acte Artistique, Université de Strasbourg;

\noindent   and Centre for Interdisciplinary Mathematical Sciences, Institute of Science, Banaras Hindu University, Vara\-nasi, 221005, India.

\noindent email:  papadop@math.unistra.fr

 \printindex

\end{document}